# Correspondences between projective planes

## June Huh

### ABSTRACT

We characterize integral homology classes of the product of two projective planes which are representable by a subvariety.

## 1. Introduction

A *planar correspondence* is a subvariety of the product of two projective planes. A substantial amount of work in the classical algebraic geometry has been devoted to the construction and analysis of such correspondences. As Fulton remarks in [Ful98, Chapter 16], a glance at the long encyclopedia article of Berzolari [Ber33] impresses one with the importance of correspondences in mathematics through the early part of the previous century. For a survey of results on planar correspondences in the classical period, see [Ber33, Chapter V] and [Sny28a].

Let $k$ be a fixed algebraically closed field. The aim of the present paper is to generalize De Jonquières' construction of planar correspondences over $k$ with controlled multidegree. This leads to a characterization of integral homology classes of $\mathbb{P}^2 \times \mathbb{P}^2$ which are representable by a (reduced and irreducible) subvariety.

THEOREM 1. *Let $\xi$ be an element in the Chow homology group*

$$\xi = a[\mathbb{P}^2 \times \mathbb{P}^0] + b[\mathbb{P}^1 \times \mathbb{P}^1] + c[\mathbb{P}^0 \times \mathbb{P}^2] \in A_2(\mathbb{P}^2 \times \mathbb{P}^2).$$

*Then $\xi$ is the class of a reduced and irreducible subvariety if and only if $a, b, c$ are nonnegative and one of the following conditions is satisfied:*

*(i)* $b > 0$ *and* $b^2 \geqslant ac$.

*(ii)* $a = 1$, $b = 0$, $c = 0$.

*(iii)* $a = 0$, $b = 0$, $c = 1$.

To the author's knowledge, the above characterization was not obtained in the classical period. For a discussion of the multidegree $a, b, c$ in the language of classical geometry, see [Bak33, SR49].

The necessity of the above conditions for the representability of $\xi$ is linked to several important achievements of the classical algebraic geometry:

(1) If $\xi$ is representable by a subvariety, then $a, b, c$ are nonnegative.

This follows from Bertini's theorem, which says that, for example, the product of two general lines transversely intersects the variety representing $\xi$ in $b$ distinct points.

(2) If $\xi$ is representable by a subvariety and $a = 1$, then $b^2 \geqslant c$.

*Keywords:* Correspondence, Projective plane, Linear system, Waring's problem.
The author was partially supported by NSF grant DMS-0943832.

For the graph of a rational map between the projective planes, the inequality $b^2 \geqslant ac$ is Bézout's theorem that two plane curves of degree $b$ without common components intersect in at most $b^2$ distinct points.

(3) If $\xi$ is representable by a subvariety, then $b^2 \geqslant ac$.

This is Hodge's index theorem on base-point-free linear systems on an algebraic surface. To be more precise, consider a resolution of singularities of the surface representing $\xi$. If $D_1, D_2$ are pull-backs of general lines from the first and the second projective plane, then the index theorem says that

$$(D_1 \cdot D_2)^2 \geqslant (D_1 \cdot D_1)(D_2 \cdot D_2).$$

(4) If $\xi$ is representable by a subvariety and $b = 0$, then $a = 1, c = 0$ or $a = 0, c = 1$.

Suppose $S$ is a subvariety representing $\xi$. By the previous item, either $c = 0$ or $a = 0$. In the former case, the image of the projection from $S$ to the second projective plane is disjoint from a line. Since this image is a projective variety contained in $\mathbb{A}^2$, the projection is constant, and hence $a = 1$.

The construction of planar correspondences with given $a, b, c$ is a classical topic. In fact, one of the first papers on the subject settles the question of existence when $a = 1, c = 1$. In [DeJ64] De Jonquières constructs a transformation of $\mathbb{P}^2$ which now bears his name, a birational transformation defined by the ratio of three homogeneous polynomials with given degree $b \geqslant 1$ and without common factors.

The present paper is devoted to the construction of planar correspondences with given $a, b, c$ which satisfy $b > 0$ and $b^2 \geqslant ac$. The possibility of such a construction is closely related to a problem on integral quadratic forms considered by Erdős, Ko, and Mordell. We briefly explain the connection.

Let us believe that a uniform construction of planar correspondences with given homology class exists. Then we should study embeddings of rational surfaces in $\mathbb{P}^2 \times \mathbb{P}^2$, because homology classes with $a = 1$ are representable only by a rational surface if they are representable by a subvariety at all. Moreover, a desingularization of any such rational surface should be a projective plane blown up at finitely many (possibly infinitely near) points.

Let $n$ be a nonnegative integer, and consider a projective plane blown up at $n$ points. We may identify its middle homology group with $\mathbb{Z}^{n+1}$, and the intersection product with

$$\mathbf{x} \circ \mathbf{y} := x_0 y_0 - x_1 y_1 - \cdots - x_n y_n, \qquad \mathbf{x}, \mathbf{y} \in \mathbb{Z}^{n+1}.$$

If $\mathbf{x} \circ \mathbf{x}$ and $\mathbf{y} \circ \mathbf{y}$ are positive, then by the reversed Cauchy-Schwarz inequality

$$(\mathbf{x} \circ \mathbf{y})^2 \geqslant (\mathbf{x} \circ \mathbf{x})(\mathbf{y} \circ \mathbf{y}).$$

In order to prove the existence of a rational surface in $\mathbb{P}^2 \times \mathbb{P}^2$ with given homology class, at the very least we should be able to answer the following question:

> Given positive integers $a, b, c$ which satisfy $b^2 \geqslant ac$, do there exist $\mathbf{x}$, $\mathbf{y}$ in $\mathbb{Z}^{n+1}$ such that $\mathbf{x} \circ \mathbf{x} = a$, $\mathbf{x} \circ \mathbf{y} = b$, $\mathbf{y} \circ \mathbf{y} = c$?

This is a subtle arithmetic question of Erdős, whose answer depends on $n$. In [Ko37] Ko answered in the negative for $n = 3$, and in the affirmative for $n = 4$. The main step in the proof of Theorem 1 is to give an affirmative answer to the same question using only those $\mathbf{x}$ and $\mathbf{y}$ which correspond to base-point-free divisors with sufficiently many sections on the blown up projective plane.

In the following section we show that the arithmetic problem above can be linearized in a suitable sense. De Jonquières' construction will be reviewed and generalized in Section 3. This generalization gives an affirmative answer to the linearized version of the arithmetic problem. In Section 4, we combine results from the previous sections to give a proof of the main theorem. A conjectural description of representable homology classes in $\mathbb{P}^3 \times \mathbb{P}^3$ will be given in Section 5.

## Acknowledgements

The author is grateful to Tsit Yuen Lam, David Leep, Bruce Reznick, and David Speyer for useful comments on Waring's problem for integral quadratic forms. He thanks Igor Dolgachev, Allen Knutson, Mircea Mustaţă, Sam Payne, and Hal Schenck for helpful discussions.

## 2. Waring's problem for integral quadratic forms

Let $n$ be a nonnegative integer, and consider the abelian group of lattice points

$$\mathbb{Z}^{n+1} := \Big\{ \mathbf{x} = (x_0, x_1, \ldots, x_n) \mid x_i \in \mathbb{Z} \Big\}.$$

Denote the Euclidean and the Lorentzian inner product on $\mathbb{Z}^{n+1}$ by

$$\mathbf{x} * \mathbf{y} := x_0 y_0 + x_1 y_1 + \cdots + x_n y_n, \qquad \mathbf{x} \circ \mathbf{y} := x_0 y_0 - x_1 y_1 - \cdots - x_n y_n.$$

By the (reversed) Cauchy-Schwarz inequality, we have

$$(\mathbf{x} * \mathbf{y})^2 \leqslant (\mathbf{x} * \mathbf{x})(\mathbf{y} * \mathbf{y}),$$

and if $\mathbf{x}$ and $\mathbf{y}$ are time-like vectors (that is, if $\mathbf{x} \circ \mathbf{x}$ and $\mathbf{y} \circ \mathbf{y}$ are positive),

$$(\mathbf{x} \circ \mathbf{y})^2 \geqslant (\mathbf{x} \circ \mathbf{x})(\mathbf{y} \circ \mathbf{y}).$$

DEFINITION 2. Let $\mathscr{L}$ be a subset of $\mathbb{Z}^{n+1}$, and let $a, b, c$ be integers.

(i) $(\mathscr{L}, *)$ *represents* $(a, b, c)$ if there exist $\mathbf{x}, \mathbf{y} \in \mathscr{L}$ such that

$$\mathbf{x} * \mathbf{x} = a, \qquad \mathbf{x} * \mathbf{y} = b, \qquad \mathbf{y} * \mathbf{y} = c.$$

$(\mathscr{L}, *)$ is *complete* if it represents every positive $(a, b, c)$ which satisfy $b^2 \leqslant ac$.

(ii) $(\mathscr{L}, \circ)$ represents $(a, b, c)$ if there exist $\mathbf{x}, \mathbf{y} \in \mathscr{L}$ such that

$$\mathbf{x} \circ \mathbf{x} = a, \qquad \mathbf{x} \circ \mathbf{y} = b, \qquad \mathbf{y} \circ \mathbf{y} = c.$$

$(\mathscr{L}, \circ)$ is *complete* if it represents every positive $(a, b, c)$ which satisfy $b^2 \geqslant ac$.

The problem of deciding the completeness of $\mathscr{L}$ may be viewed as an extension of Waring's problem [Ko37, Mor30, Mor32]. For example, $(\mathscr{L}, *)$ represents $(a, b, c)$ if and only if every binary quadratic form

$$ax^2 + 2bxy + cy^2$$

is a sum of $k$ squares of integral linear forms with coefficients from $\mathscr{L}$.

*Example* 3. $(\mathbb{Z}^4, *)$ is not complete because the condition fails for $a = 1, b = 2, c = 19$. To see this, note that 19 is a sum of four squares in exactly two different ways

$$19 = 4^2 + 1^2 + 1^2 + 1^2 = 3^2 + 3^2 + 1^2 + 0^2.$$

Mordell proved in [Mor30, Mor32] that $(\mathbb{Z}^5, *)$ is complete.

*Example* 4. $(\mathbb{N}^8, *)$ is not complete because the condition fails for $a = 8, b = 1, c = 8$. To see this, note that 8 is a sum of eight squares in exactly three different ways

$$8 = 2^2 + 2^2 + 0^2 + 0^2 + 0^2 + 0^2 + 0^2 + 0^2$$
$$= 2^2 + 1^2 + 1^2 + 1^2 + 1^2 + 0^2 + 0^2 + 0^2$$
$$= 1^2 + 1^2 + 1^2 + 1^2 + 1^2 + 1^2 + 1^2 + 1^2.$$

The author does not know the smallest $n$ for which $(\mathbb{N}^{n+1}, *)$ is complete.

Showing that $(\mathscr{L}, \circ)$ is or is not complete is a delicate problem in general. Ko proved in [Ko37] that $(\mathbb{Z}^4, \circ)$ is not complete and that $(\mathbb{Z}^5, \circ)$ is complete, thus answering a question of Erdős.

Proposition 6 below shows that the problem of representation by $(\mathscr{L}, \circ)$ can be linearized in a suitable sense if $\mathscr{L}$ is closed under the addition of $\mathbb{Z}^{n+1}$. This observation will play an important role in the proof of Theorem 1.

DEFINITION 5. Let $k$ be a positive integer. $(\mathscr{L}, \circ)$ is *linearly $k$-complete* if it represents every positive $(a, b, c)$ which satisfy

$$2b \geqslant a + c \quad \text{and} \quad b \leqslant k.$$

Similarly, we say that $(\mathscr{L}, \circ)$ is *$k$-complete* if it represents every positive $(a, b, c)$ which satisfy

$$b^2 \geqslant ac \quad \text{and} \quad b \leqslant k.$$

By the inequality of arithmetic and geometric means, if $(\mathscr{L}, \circ)$ is $k$-complete, then $(\mathscr{L}, \circ)$ is linearly $k$-complete.

PROPOSITION 6. *If $(\mathscr{L}, \circ)$ is linearly $k$-complete and $\mathscr{L}$ is closed under the addition of $\mathbb{Z}^{n+1}$, then $(\mathscr{L}, \circ)$ is $k$-complete.*

*Proof.* Let $a, b, c$ be positive integers which satisfy $b^2 \geqslant ac$ and $b \leqslant k$. We show by induction on $k$ that there exist $\mathbf{x}, \mathbf{y} \in \mathscr{L}$ such that

$$\mathbf{x} \circ \mathbf{x} = a, \qquad \mathbf{x} \circ \mathbf{y} = b, \qquad \mathbf{y} \circ \mathbf{y} = c.$$

The base case $k = 1$ is immediate.

We may suppose that $2b < a + c$ and $c \geqslant a$. Since $b^2 \geqslant ac$, this implies that $b \geqslant a$. The idea is to rewrite $a, b, c$ by

$$a = a, \qquad b = a + (b - a), \qquad c = a + 2(b - a) + (a + c - 2b),$$

and consider the new triple

$$\tilde{a} := a, \qquad \tilde{b} := b - a, \qquad \tilde{c} := a + c - 2b.$$

Under our assumptions on $a, b, c$, the new triple $\tilde{a}, \tilde{b}, \tilde{c}$ has the following properties:

(1) $\tilde{a}$, $\tilde{b}$, $\tilde{c}$ are positive integers.
(2) The discriminant remains unchanged: $\tilde{b}^2 - \tilde{a}\tilde{c} = b^2 - ac \geqslant 0$.
(3) The induction invariant drops: $\tilde{b} < b \leqslant k$.

Therefore we may use the induction hypothesis to find $\tilde{\mathbf{x}}, \tilde{\mathbf{y}} \in \mathscr{L}$ such that

$$\tilde{\mathbf{x}} \circ \tilde{\mathbf{x}} = \tilde{a}, \qquad \tilde{\mathbf{x}} \circ \tilde{\mathbf{y}} = \tilde{b}, \qquad \tilde{\mathbf{y}} \circ \tilde{\mathbf{y}} = \tilde{c}.$$

Since $\circ$ is bilinear, we have

$$\tilde{\mathbf{x}} \circ \tilde{\mathbf{x}} = a, \qquad \tilde{\mathbf{x}} \circ (\tilde{\mathbf{x}} + \tilde{\mathbf{y}}) = b, \qquad (\tilde{\mathbf{x}} + \tilde{\mathbf{y}}) \circ (\tilde{\mathbf{x}} + \tilde{\mathbf{y}}) = c.$$

Now $\mathbf{x} := \tilde{\mathbf{x}}$ and $\mathbf{y} := \tilde{\mathbf{x}} + \tilde{\mathbf{y}}$ are elements of $\mathscr{L}$ with the desired properties. $\qquad\square$

As an application, we show that $(\mathbb{N}^7, \circ)$ is complete. A modification of the argument below will play a role in the proof of Theorem 1.

COROLLARY 7. $(\mathbb{N}^7, \circ)$ *is complete.*

*Proof.* By Proposition 6, it is enough to show that $(\mathbb{N}^7, \circ)$ is linearly $k$-complete for all $k$. Let $a, b, c$ be positive integers which satisfy $2b \geqslant a + c$. We may suppose that $c \leqslant a$ and $c \leqslant b$.

Define nonnegative integers

$$r_1 := \lfloor \frac{c}{2} \rfloor, \qquad r_2 := b - c, \qquad r_3 := 2b - a - c.$$

Use Lagrange's four squares theorem to find nonnegative integers $n_1, n_2, n_3, n_4$ such that

$$r_3 = n_1^2 + n_2^2 + n_3^2 + n_4^2.$$

If $c$ is odd, then set

$$\mathbf{x} := (r_1 + r_2 + 1, r_1 + r_2, 0, n_1, n_2, n_3, n_4), \qquad \mathbf{y} := (r_1 + 1, r_1, 0, 0, 0, 0, 0).$$

We have

$$\mathbf{x} \circ \mathbf{x} = (2r_1 + 1) + 2r_2 - r_3, \quad \mathbf{x} \circ \mathbf{y} = (2r_1 + 1) + r_2, \quad \mathbf{y} \circ \mathbf{y} = 2r_1 + 1.$$

If $c$ is even, then set

$$\mathbf{x} := (r_1 + r_2 + 1, r_1 + r_2, 1, n_1, n_2, n_3, n_4), \qquad \mathbf{y} := (r_1 + 1, r_1, 1, 0, 0, 0, 0).$$

We have

$$\mathbf{x} \circ \mathbf{x} = 2r_1 + 2r_2 - r_3, \quad \mathbf{x} \circ \mathbf{y} = 2r_1 + r_2, \quad \mathbf{y} \circ \mathbf{y} = 2r_1.$$

In both cases,

$$\mathbf{x} \circ \mathbf{x} = a, \qquad \mathbf{x} \circ \mathbf{y} = b, \qquad \mathbf{y} \circ \mathbf{y} = c.$$

$\qquad\square$

The author does not know the smallest $n$ for which $(\mathbb{N}^{n+1}, \circ)$ is complete.

## 3. Linear systems of De Jonquières type

Let $\mathbf{p} = (p_1, p_2, \ldots, p_n)$ be a sequence of distinct points in the projective plane $\mathbb{P}^2$. Consider the set of nonnegative lattice points

$$\mathbb{N}^{n+1} := \Big\{ \mathbf{m} := (d, m_1, m_2, \ldots, m_n) \mid d \geqslant 0, \ m_i \geqslant 0 \Big\}.$$

An element $\mathbf{m} \in \mathbb{N}^{n+1}$ as above, together with the sequence $\mathbf{p}$, defines a linear system of plane curves

$$L(\mathbf{p}, \mathbf{m}) := \Big\{ C \mid \deg C = d \text{ and the multiplicity of } C \text{ at } p_i \text{ is at least } m_i \text{ for all } i \Big\}.$$

Note that every linear system of $\mathbb{P}^2$ is of the form $L(\mathbf{p}, \mathbf{m})$ for some $\mathbf{p}$, $\mathbf{m}$, and $n$.

DEFINITION 8. $L(\mathbf{p}, \mathbf{m})$ has *no unassigned base points* if

(i) $L(\mathbf{p}, \mathbf{m})$ is nonempty,

(ii) no point other than $p_1, p_2, \ldots, p_n$ is contained in $C$ for all $C \in L(\mathbf{p}, \mathbf{m})$,

(iii) no line is contained in the tangent cone of $C$ at $p_i$ for all $C \in L(\mathbf{p}, \mathbf{m})$, and

(iv) there is an element of $L(\mathbf{p}, \mathbf{m})$ which has multiplicity $m_i$ at $p_i$ for all $i$.

In other words, we require that the linear system has no base points other than $p_1, \ldots, p_n$, both proper or infinitely near, and a general member of the linear system has the expected multiplicity at $p_i$ for all $i$.

When nonempty, we may view $L(\mathbf{p}, \mathbf{m})$ as a projective space. We denote the rational map associated to $L(\mathbf{p}, \mathbf{m})$ by

$$\varphi(\mathbf{p}, \mathbf{m}) : \mathbb{P}^2 \longrightarrow L(\mathbf{p}, \mathbf{m})^{\vee}.$$

If $p$ is not a base point of the linear system, then $\varphi(\mathbf{p}, \mathbf{m})$ maps $p$ to the hyperplane of curves passing through $p$.

DEFINITION 9. $L(\mathbf{p}, \mathbf{m})$ is *very big* if $\varphi(\mathbf{p}, \mathbf{m})$ maps its domain birationally onto its image.

The notion is a birational analogue of *very ample*. The following is the main result of this section.

PROPOSITION 10. *Define*

$$\mathscr{L}(\mathbf{p}) := \Big\{ \mathbf{m} \in \mathbb{N}^{n+1} \mid L(\mathbf{p}, \mathbf{m}) \text{ is very big and has no unassigned base points} \Big\}.$$

*Then* $\big(\mathscr{L}(\mathbf{p}), \circ\big)$ *is $k$-complete for $k = \lfloor n/2 \rfloor$ and a sufficiently general $\mathbf{p}$.*

The proof of Proposition 10 is built upon results of De Jonquières [DeJ64, DeJ85]. We recall the construction and the needed properties of De Jonquières transformation of $\mathbb{P}^2$. Modern treatments can be found in [Dol12, Chapter 7] and [KSC04, Chapter 2].

A *De Jonquières transformation* of degree $d \geqslant 1$ is a birational map of the form $\varphi(\mathbf{p}, \mathbf{m})$, where

$$\mathbf{m} = (d, d-1, \underbrace{1, 1, \ldots, 1}_{2d-2}, \underbrace{0, 0, \ldots, 0}_{n-2d+1}).$$

The result of De Jonquières is that, for $\mathbf{m}$ as above and a sufficiently general $\mathbf{p}$,

(1) $\dim L(\mathbf{p}, \mathbf{m}) = 2$,

(2) $L(\mathbf{p}, \mathbf{m})$ has no unassigned base points, and

(3) $\varphi(\mathbf{p}, \mathbf{m})$ is a birational transformation of $\mathbb{P}^2$.

It is necessary to assume that $\mathbf{p}$ is sufficiently general. For example, if $d = 3$ and $p_2, p_3, p_4, p_5$ are collinear, then all three conditions above fail to hold for $L(\mathbf{p}, \mathbf{m})$.

*Remark* 11. Interesting De Jonquières transformations can be obtained by allowing some of the base points to be infinitely near. We will not need this extension.

Lemma 12 and Lemma 13 below will be needed in the proof of Proposition 10.

LEMMA 12. *Define* $\mathbf{m} = (d, m_1, m_2, \ldots, m_n) \in \mathbb{N}^{n+1}$ *to be of De Jonquières type if*

(i) $d \geqslant 1$, $n \geqslant 2d - 1$, *and* $m_1 = d - 1$,

(ii) $m_2, m_3, \ldots, m_n$ *are either zero or one, and*

(iii) *at most* $2d - 2$ *among* $m_2, m_3, \ldots, m_n$ *are nonzero.*

If $\mathbf{m}$ is of De Jonquières type, then $L(\mathbf{p}, \mathbf{m})$ is very big and has no unassigned base points for a sufficiently general $\mathbf{p}$.

*Proof.* Let $r$ be the number of nonzeros among $m_2, m_3, \ldots, m_n$. We may suppose that $n = 2d-1$ and

$$\mathbf{m} = (d, d-1, \underbrace{1, 1, \ldots, 1}_{r}, \underbrace{0, 0, \ldots, 0}_{2d-2-r}).$$

Define

$$\mathbf{n} := (d, d-1, \underbrace{1, 1, \ldots, 1, 1, 1, \ldots, 1}_{2d-2}).$$

There is an inclusion between the linear systems

$$\iota : L(\mathbf{p}, \mathbf{n}) \longrightarrow L(\mathbf{p}, \mathbf{m}).$$

For any $\mathbf{p}$, we have the commutative diagram of rational maps

$$\begin{array}{ccc} \mathbb{P}^2 & \xrightarrow{\varphi(\mathbf{p}, \mathbf{m})} & L(\mathbf{p}, \mathbf{m})^\vee \\ & \searrow^{\varphi(\mathbf{p}, \mathbf{n})} & \downarrow^{\iota^\vee} \\ & & L(\mathbf{p}, \mathbf{n})^\vee. \end{array}$$

By the result of De Jonquières, we may choose $\mathbf{p}$ sufficiently general so that $\varphi(\mathbf{p}, \mathbf{n})$ is a birational transformation of $\mathbb{P}^2$. Then the commutative diagram shows that $L(\mathbf{p}, \mathbf{m})$ is very big.

Next we show that $L(\mathbf{p}, \mathbf{m})$ has no unassigned base points. Let $p$ be a point different from $p_1, \ldots, p_n$. Again by the result of De Jonquières, there is a sequence of distinct points $\mathbf{q} = (q_1, \ldots, q_n)$ such that

(1) $\varphi(\mathbf{q}, \mathbf{n})$ is a De Jonquières transformation of $\mathbb{P}^2$,

(2) $q_i$ is equal to $p_i$ for $1 \leqslant i \leqslant r+1$, and

(3) $q_i$ is different from $p$ and from $p_{r+1}, \ldots, p_n$ for $r+1 < i \leqslant n$.

Note that $L(\mathbf{q}, \mathbf{n})$ is a subspace of $L(\mathbf{p}, \mathbf{m})$ which has no unassigned base points. Since $p$ can be any point different from $p_1, \ldots, p_n$, it follows that $L(\mathbf{p}, \mathbf{m})$ has no unassigned base points. $\qquad\square$

LEMMA 13. *Let $\mathbf{m}_1, \mathbf{m}_2 \in \mathbb{N}^{n+1}$ and $\mathbf{m}_3 = \mathbf{m}_1 + \mathbf{m}_2$.*

(i) *If $L(\mathbf{p}, \mathbf{m}_1)$ is very big and $L(\mathbf{p}, \mathbf{m}_2)$ is nonempty, then $L(\mathbf{p}, \mathbf{m}_3)$ is very big.*

(ii) *If $L(\mathbf{p}, \mathbf{m}_1)$ and $L(\mathbf{p}, \mathbf{m}_2)$ have no unassigned base points, then $L(\mathbf{p}, \mathbf{m}_3)$ has no unassigned base points.*

*Proof.* Choose an element $C_2$ of the linear system $L(\mathbf{p}, \mathbf{m}_2)$. This defines an embedding

$$\iota_2 : L(\mathbf{p}, \mathbf{m}_1) \longrightarrow L(\mathbf{p}, \mathbf{m}_3), \qquad C_1 \longmapsto C_1 \cup C_2$$

and the commutative diagram of rational maps

$$
\begin{array}{ccc}
\mathbb{P}^2 & \xrightarrow{\;\varphi(\mathbf{p},\mathbf{m}_3)\;} & L(\mathbf{p},\mathbf{m}_3)^\vee \\
& \varphi(\mathbf{p},\mathbf{m}_1)\searrow & \downarrow \iota_2^\vee \\
& & L(\mathbf{p},\mathbf{m}_1)^\vee .
\end{array}
$$

Since $L(\mathbf{p},\mathbf{m}_1)$ is very big, the commutative diagram shows that $L(\mathbf{p},\mathbf{m}_3)$ is very big.

Next we show that $L(\mathbf{p},\mathbf{m}_3)$ has no unassigned base points. Consider the subset

$$
L(\mathbf{p},\mathbf{m}_1) + L(\mathbf{p},\mathbf{m}_2) := \Big\{ C_1 \cup C_2 \mid C_1 \in L(\mathbf{p},\mathbf{m}_1),\ C_2 \in L(\mathbf{p},\mathbf{m}_2) \Big\} \subseteq L(\mathbf{p},\mathbf{m}_3).
$$

Since $L(\mathbf{p},\mathbf{m}_1)$ and $L(\mathbf{p},\mathbf{m}_2)$ have no unassigned base points, there is an element of the above subset which has the expected multiplicity at $p_1,\ldots,p_n$ and which does not contain a given point different from $p_1,\ldots,p_n$, proper or infinitely near. It follows that $L(\mathbf{p},\mathbf{m}_3)$ has no unassigned base points. $\qquad\square$

*Proof of Proposition 10.* By Lemma 13, we know that $\mathscr{L}(\mathbf{p})$ is closed under the addition of $\mathbb{Z}^{n+1}$ for any $\mathbf{p}$. Therefore, by Proposition 6, it is enough to prove that $\mathscr{L}(\mathbf{p})$ is linearly $k$-complete for a sufficiently general $\mathbf{p}$. This linear version of the problem can be solved directly by using linear systems of De Jonquières type.

Recall that $\mathbf{m} = (d,m_1,m_2,\ldots,m_n) \in \mathbb{N}^{n+1}$ is said to be *of De Jonquières type* if

(i) $d \geqslant 1$, $n \geqslant 2d-1$, and $m_1 = d-1$,

(ii) $m_2,m_3,\ldots,m_n$ are either zero or one, and

(iii) at most $2d-2$ among $m_2,m_3,\ldots,m_n$ are nonzero.

Let $\mathscr{D}$ be the set of all elements of De Jonquières type in $\mathbb{N}^{n+1}$. Since $n \geqslant 2d-1$, $\mathscr{D}$ has only finitely many elements. Therefore, by Lemma 12, we may choose $\mathbf{p}$ sufficiently general so that

$$
\mathscr{D} \subseteq \mathscr{L}(\mathbf{p}).
$$

Let $a,b,c$ be positive integers which satisfy $2b \geqslant a+c$ and $b \leqslant k$. We may suppose that $c \leqslant a$ and $c \leqslant b$. We show that there exist $\mathbf{m}_1,\mathbf{m}_2 \in \mathscr{D}$ which satisfy

$$
\mathbf{m}_1 \circ \mathbf{m}_1 = a, \qquad \mathbf{m}_1 \circ \mathbf{m}_2 = b, \qquad \mathbf{m}_2 \circ \mathbf{m}_2 = c.
$$

For this we mimic the proof of Corollary 7. Define nonnegative integers

$$
r_1 := \lfloor \tfrac{c}{2} \rfloor, \qquad r_2 := b-c, \qquad r_3 := 2b-a-c.
$$

Note that $r_2 \leqslant n-2$.

(1) If $c$ is odd, then set

$$
\mathbf{m}_1 := (r_1+r_2+1, r_1+r_2, 0, \underbrace{1,1,\ldots,1}_{r_3}, \underbrace{0,0,\ldots,0}_{n-2-r_3}),
$$

$$
\mathbf{m}_2 := (r_1+1, r_1, 0, \underbrace{0,0,\ldots,0}_{r_3}, \underbrace{0,0,\ldots,0}_{n-2-r_3}).
$$

It is easy to check that

$$
r_3 \leqslant 2r_1 + 2r_2 < n.
$$

The two inequalities show that $\mathbf{m}_1, \mathbf{m}_2$ are of De Jonquières type. We have
$$\mathbf{m}_1 \circ \mathbf{m}_1 = (2r_1 + 1) + 2r_2 - r_3, \quad \mathbf{m}_1 \circ \mathbf{m}_2 = (2r_1 + 1) + r_2, \quad \mathbf{m}_2 \circ \mathbf{m}_2 = 2r_1 + 1.$$

(2) If $c$ is even, then set
$$\mathbf{m}_1 := (r_1 + r_2 + 1, r_1 + r_2, 1, \underbrace{1, 1, \ldots, 1}_{r_3}, \underbrace{0, 0, \ldots, 0}_{n-2-r_3}),$$
$$\mathbf{m}_2 := (r_1 + 1, r_1, 1, \underbrace{0, 0, \ldots, 0}_{r_3}, \underbrace{0, 0, \ldots, 0}_{n-2-r_3}).$$

Since $c$ is even,
$$r_3 + 1 \leqslant 2r_1 + 2r_2 < n \quad \text{and} \quad 1 \leqslant 2r_1.$$

The three inequalities show that $\mathbf{m}_1, \mathbf{m}_2$ are of De Jonquières type. We have
$$\mathbf{m}_1 \circ \mathbf{m}_1 = 2r_1 + 2r_2 - r_3, \quad \mathbf{m}_1 \circ \mathbf{m}_2 = 2r_1 + r_2, \quad \mathbf{m}_2 \circ \mathbf{m}_2 = 2r_1.$$

In both cases,
$$\mathbf{m}_1 \circ \mathbf{m}_1 = a, \qquad \mathbf{m}_1 \circ \mathbf{m}_2 = b, \qquad \mathbf{m}_2 \circ \mathbf{m}_2 = c.$$

$\square$

## 4. Proof of Theorem 1

We first characterize homology classes of reduced and irreducible surfaces in $\mathbb{P}^2 \times \mathbb{P}^1$.

Lemma 14. *Let $\xi$ be an element in the Chow homology group*
$$\xi = a[\mathbb{P}^2 \times \mathbb{P}^0] + b[\mathbb{P}^1 \times \mathbb{P}^1] \in A_2(\mathbb{P}^2 \times \mathbb{P}^1).$$

*Then $\xi$ is the class of a reduced and irreducible subvariety if and only if one of the following conditions is satisfied:*

(i) $b > 0$ and $a \geqslant 0$.

(ii) $b = 0$ and $a = 1$.

*Proof.* Suppose $\xi$ is the class of a reduced and irreducible subvariety $X$. Then $X$ is defined by an irreducible bihomogeneous polynomial in two sets of variables $z_0, z_1, z_2$ and $w_0, w_1$ with respective degrees $b$ and $a$. The assertion that $b = 0$ implies $a = 1$ is precisely the fundamental theorem of algebra applied to the defining equation of $X$.

Conversely, if $b > 0$ and $a \geqslant 0$, then a sufficiently general bihomogeneous polynomial in variables $z_0, z_1, z_2$ and $w_0, w_1$ with respective degrees $b$ and $a$ is irreducible by Bertini's theorem. This proves Lemma 14. $\square$

Representable homology classes in $A_2(\mathbb{P}^1 \times \mathbb{P}^2)$ can be characterized in the same way.

*Remark* 15. Let $\xi$ be an element in the Chow homology group
$$\xi = a[\mathbb{P}^1 \times \mathbb{P}^0] + b[\mathbb{P}^0 \times \mathbb{P}^1] \in A_1(\mathbb{P}^1 \times \mathbb{P}^1).$$

Then $\xi$ is the class of a reduced and irreducible subvariety if and only if one of the following conditions is satisfied:

(i) $a > 0$ and $b > 0$.

(ii) $a = 1$ and $b = 0$.

(iii) $a = 0$ and $b = 1$.

The proof is similar to that of Lemma 14.

Let $X$ be a reduced and irreducible surface, and let $f : X \longrightarrow \mathbb{P}^{N_1} \times \mathbb{P}^{N_2}$ be a regular map to a biprojective space with $N_1 \geqslant 2$, $N_2 \geqslant 2$. We denote the two projections by

$$
\begin{array}{ccc}
 & X & \\
\text{pr}_1 \swarrow & & \searrow \text{pr}_2 \\
\mathbb{P}^{N_1} & & \mathbb{P}^{N_2}.
\end{array}
$$

The following result will serve as a final preparation for the proof of Theorem 1.

Lemma 16. *Consider the commutative diagram of rational maps*

$$
\begin{array}{ccc}
X & \xrightarrow{\;f\;} & \mathbb{P}^{N_1} \times \mathbb{P}^{N_2} \\
 & \pi \searrow & \downarrow \pi_1 \times \pi_2 \\
 & & \mathbb{P}^2 \times \mathbb{P}^2
\end{array}
$$

*where $\pi_1$ and $\pi_2$ are independently chosen general linear projections. Then $\pi$ is a regular map, and if $\text{pr}_1$ and $\text{pr}_2$ map $X$ birationally onto $\text{pr}_1(X)$ and $\text{pr}_2(X)$ respectively, then $\pi$ maps $X$ birationally onto $\pi(X)$.*

It is not enough to assume that $\text{pr}_1$ maps $X$ birationally onto $\text{pr}_1(X)$. For example, if $\text{pr}_1$ is an embedding of a degree $d$ surface in $\mathbb{P}^3$ and $\text{pr}_2$ is a constant map to $\mathbb{P}^3$, then $\pi : X \longrightarrow \pi(X)$ has degree $d$ for sufficiently general $\pi_1$ and $\pi_2$.

*Proof.* The center of the linear projection $\pi_1$ (respectively $\pi_2$) is either empty or has codimension 3 in $\mathbb{P}^{N_1}$ (respectively in $\mathbb{P}^{N_2}$). Therefore $\pi$ is defined everywhere on $X$ for sufficiently general $\pi_1$ and $\pi_2$.

Suppose $\text{pr}_1$ and $\text{pr}_2$ map $X$ birationally onto $\text{pr}_1(X)$ and $\text{pr}_2(X)$ respectively. Define $\widetilde{f}$, $g$, and $h$ by the commutative diagram

$$
\begin{array}{ccc}
X & \xrightarrow{\;f\;} & \mathbb{P}^{N_1} \times \mathbb{P}^{N_2} \\
g \downarrow \quad \widetilde{f} \searrow & & \downarrow \pi_1 \times \text{id} \\
 & \mathbb{P}^2 \times \mathbb{P}^{N_2} & \\
 & h \swarrow & \searrow \\
\mathbb{P}^2 & & \mathbb{P}^{N_2}.
\end{array}
$$

We claim that

(1) $\widetilde{f}$ maps $X$ birationally onto $\widetilde{f}(X)$, and

(2) $\widetilde{h} := h|_{\widetilde{f}(X)}$ has a reduced general fiber for a sufficiently general $\pi_1$.

The first assertion is valid for any $\pi_1$ because $\widetilde{f}$ is a factor of $\text{pr}_2$. For the second assertion, note that a general codimension 2 linear subspace of $\mathbb{P}^{N_1}$ intersects $\text{pr}_1(X)$ in finitely many reduced points. Since $X$ is mapped birationally onto $\text{pr}_1(X)$, the previous sentence implies that $g$ has a reduced general fiber for a sufficiently general $\pi_1$. Therefore $\widetilde{h}$ has a reduced general fiber for a sufficiently general $\pi_1$.

We show that $\pi$ maps $X$ birationally onto $\pi(X)$ by induction on $N_2$. Suppose $N_2 > 2$, and let $\widetilde{f}$ and $\widetilde{h}$ be as above. Consider the linear projection $\mathrm{p}_2 : \mathbb{P}^{N_2} \longrightarrow \mathbb{P}^{N_2-1}$ from a point $y$, and define $i$ by the commutative diagram

$$
\begin{array}{ccc}
\widetilde{f}(X) & \longrightarrow & \mathbb{P}^2 \times \mathbb{P}^{N_2} \\
& i \searrow & \downarrow \mathrm{id}\times\mathrm{p}_2 \\
\widetilde{h} \downarrow & & \mathbb{P}^2 \times \mathbb{P}^{N_2-1} \\
& \swarrow & \\
\mathbb{P}^2. & &
\end{array}
$$

If $x_1$ is a sufficiently general point of $\mathbb{P}^2$, then the fiber of $x_1$ over $\widetilde{h}$ is a reduced set of points

$$
\widetilde{h}^{-1}(x_1) = \Big\{ (x_1, y_1), \ldots, (x_1, y_m) \Big\}, \qquad y_1, \ldots, y_m \in \mathbb{P}^{N_2}.
$$

If the center $y$ is not contained in the union of the lines joining $y_i$ and $y_j$, and if $\mathbb{P}^2 \times \{y\}$ is disjoint from $\widetilde{f}(X)$, then

$$
i^{-1}\Big( i(x_1, y_1) \Big) = \Big\{ (x_1, y_1) \Big\}.
$$

It follows that $i$ maps $\widetilde{f}(X)$ birationally onto its image in $\mathbb{P}^2 \times \mathbb{P}^{N_2-1}$ for a sufficiently general $y$. The proof is completed by induction. $\qquad\square$

*Proof of Theorem 1.* We construct a reduced and irreducible surface in $\mathbb{P}^2 \times \mathbb{P}^2$ with given $a, b, c$ which satisfy $b > 0$ and $b^2 \geqslant ac$. If $a = 0$ or $c = 0$, then we may use Lemma 14.

Suppose $a, b, c$ are positive. Let $X$ be the blowup of $\mathbb{P}^2$ at $n \geqslant 2b$ sufficiently general points. By Proposition 10, there are base-point-free divisors $D_1, D_2$ of $X$ such that

$$
D_1 \cdot D_1 = a, \qquad D_1 \cdot D_2 = b, \qquad D_2 \cdot D_2 = c,
$$

whose linear systems map $X$ birationally onto their respective images. Let $L_1, L_2$ be the linear systems of $D_1, D_2$, and write $\varphi_1, \varphi_2$ for the corresponding rational maps. We apply Lemma 16 to the product

$$
\varphi_1 \times \varphi_2 : X \longrightarrow L_1^\vee \times L_2^\vee \simeq \mathbb{P}^{N_1} \times \mathbb{P}^{N_2}.
$$

If $\widetilde{L}_1, \widetilde{L}_2$ are sufficiently general two dimensional linear subspaces of $L_1, L_2$ respectively, then the biprojection

$$
\widetilde{\varphi}_1 \times \widetilde{\varphi}_2 : X \longrightarrow \widetilde{L}_1^\vee \times \widetilde{L}_2^\vee \simeq \mathbb{P}^2 \times \mathbb{P}^2
$$

is a regular map which maps $X$ birationally onto its image. By the projection formula [Ful98, Example 2.4.3], we have

$$
\big[ (\widetilde{\varphi}_1 \times \widetilde{\varphi}_2)(X) \big] = a[\mathbb{P}^2 \times \mathbb{P}^0] + b[\mathbb{P}^1 \times \mathbb{P}^1] + c[\mathbb{P}^0 \times \mathbb{P}^2] \in A_2(\mathbb{P}^2 \times \mathbb{P}^2).
$$

$\qquad\square$

*Remark* 17. The proof of Theorem 1 shows that a (reduced and irreducible) surface in $\mathbb{P}^2 \times \mathbb{P}^2$ is homologous to either

(i) a surface in $\mathbb{P}^1 \times \mathbb{P}^2$,

(ii) a surface in $\mathbb{P}^2 \times \mathbb{P}^1$, or

(iii) the image of a blown-up projective plane whose embedding is built up from linear systems of De Jonquières type.

It is pleasant to recall the classical fact that De Jonquières transformations are basic building blocks of birational transformations of $\mathbb{P}^2$. See, for example, [KSC04, Theorem 2.30].

*Remark* 18. Let $X$ be a smooth projective variety. In [Har74, Question 1.3] Hartshorne asks which homology classes of $X$ can be represented by an irreducible nonsingular subvariety. The author does not know whether the characterization of representability in Theorem 1 remains unchanged if one requires subvarieties to be nonsingular. We note that, when $X$ is a complex Grassmannian, there are subvarieties of $X$ which are not smoothable up to homological equivalence [Bry10, Cos11, HRT74, Hon05].

## 5. Further discussion

We conjecture that an analogue of Theorem 1 remains valid in dimension 3. For a survey of results on three dimensional correspondences in the classical period, see [Ber33, Sny28b, Sny34].

CONJECTURE 19. *Let $\xi$ be an element in the Chow homology group*

$$\xi = a[\mathbb{P}^3 \times \mathbb{P}^0] + b[\mathbb{P}^2 \times \mathbb{P}^1] + c[\mathbb{P}^1 \times \mathbb{P}^2] + d[\mathbb{P}^0 \times \mathbb{P}^3] \in A_3(\mathbb{P}^3 \times \mathbb{P}^3).$$

*Then $\xi$ is the class of a reduced and irreducible subvariety if and only if $a, b, c, d$ are nonnegative and one of the following conditions is satisfied:*

(i) $b^2 + c^2 > 0$ and $b^2 \geqslant ac$ and $c^2 \geqslant bd$.

(ii) $a = 1$, $b = 0$, $c = 0$, $d = 0$.

(iii) $a = 0$, $b = 0$, $c = 0$, $d = 1$.

The necessity of the above numerical conditions for the representability of $\xi$ follows from Theorem 20 below. For the sufficiency of the numerical conditions in the case of Cremona transformations (that is, when $a = 1$ and $d = 1$), see [Pan]. The construction is based on a 3-dimensional generalization of the De Jonquières birational transformation [Pan00, Pan01].

To illustrate the nature of Conjecture 19, we quote below [Huh12, Theorem 21] which characterizes representable homology classes in $\mathbb{P}^n \times \mathbb{P}^m$ for any nonnegative $n, m$, up to an integral multiple. Recall that a sequence $e_0, e_1, \ldots, e_n$ of integers is said to be *log-concave* if for all $0 < i < n$,

$$e_{i-1}e_{i+1} \leqslant e_i^2,$$

and it is said to have *no internal zeros* if there do not exist $i < j < k$ satisfying

$$e_i \neq 0, \quad e_j = 0, \quad e_k \neq 0.$$

THEOREM 20. *Let $\xi$ be an element in the Chow homology group*

$$\xi = \sum_i e_i[\mathbb{P}^i \times \mathbb{P}^{k-i}] \in A_k(\mathbb{P}^n \times \mathbb{P}^m).$$

(i) *If $\xi$ is an integer multiple of either*

$$[\mathbb{P}^m \times \mathbb{P}^n], [\mathbb{P}^m \times \mathbb{P}^0], [\mathbb{P}^0 \times \mathbb{P}^n], [\mathbb{P}^0 \times \mathbb{P}^0],$$

*then $\xi$ is the class of a reduced and irreducible subvariety if and only if the integer is 1.*

(ii) *If otherwise, some positive integer multiple of $\xi$ is the class of a reduced and irreducible subvariety if and only if the $e_i$ form a nonzero log-concave sequence of nonnegative integers with no internal zeros.*

It is necessary to take a positive integer multiple of $\xi$ in the second part of Theorem 20. A result of Pirio and Russo [PRa, Corollary 5.3] implies that there is no reduced and irreducible subvariety of $\mathbb{P}^5 \times \mathbb{P}^5$ which has the homology class

$$1[\mathbb{P}^5 \times \mathbb{P}^0] + 2[\mathbb{P}^4 \times \mathbb{P}^1] + 3[\mathbb{P}^3 \times \mathbb{P}^2] + 4[\mathbb{P}^2 \times \mathbb{P}^3] + 2[\mathbb{P}^1 \times \mathbb{P}^4] + 1[\mathbb{P}^0 \times \mathbb{P}^5] \in A_5(\mathbb{P}^5 \times \mathbb{P}^5).$$

Note that the sequence $(1, 2, 3, 4, 2, 1)$ is log-concave and has no internal zeros. The nonexistence can (also) be deduced from an explicit classification of quadro-quadric Cremona transformations in dimension five [PRb, Table 10].

In general, it is difficult to characterize homology classes of subvarieties of a given algebraic variety, even when the ambient variety is a smooth projective toric variety over $\mathbb{C}$. For example, when the complex toric variety is the one corresponding to the $n$-dimensional permutohedron, the problem of characterizing representable homology classes is at least as difficult as identifying matroids with $n + 1$ elements representable over the complex numbers [Fin12, KP11]. The latter is a difficult problem in a rather precise sense, and one does not expect an answer which is uniform with respect to $n$ [MNW, Vam78].

For what is known about representable homology classes in homogeneous varieties, see [Bry10, Cos11, CR, Hon05, Hon07, Per02] and references therein.

June Huh    junehuh@umich.edu

Department of Mathematics, University of Michigan,  Ann Arbor, MI 48109,  USA